\documentclass[a4paper, 12pt]{amsart}
\usepackage{amsmath, amssymb, eucal, amscd, amstext, enumerate}
\usepackage{amsfonts,amsthm}
\usepackage{mathrsfs}

\usepackage{color,hyperref}

\theoremstyle{plain}
\newtheorem{thm}{Theorem}[section]
\newtheorem{lem}[thm]{Lemma}
\newtheorem{eg}[thm]{Example}

\newtheorem{defn}[thm]{Definition}
\newtheorem{rem}[thm]{Remark}
\newtheorem{rem-ntn}[thm]{Remark and Notation}

\newenvironment{prf}{{\noindent \textbf{Proof:}\ }}{\hfill $\Box$\\ \smallskip}

\numberwithin{equation}{section}

\newcommand{\rd}{{\rm d}}

\newcommand{\abs}[1]{\left\vert#1\right\vert}
\newcommand{\norm}[1]{\left\|#1\right\|}

\newcommand{\CF}{\mathcal{F}}

\newcommand{\CA}{\mathcal{A}}

\newcommand{\CM}{\mathcal{M}}

\newcommand{\BR}{\mathbb{R}}
\newcommand{\BN}{\mathbb{N}}
\newcommand{\BK}{\mathbb{K}}
\newcommand{\BC}{\mathbb{C}}
\newcommand{\BZ}{\mathbb{Z}}

\usepackage{color}

\begin{document}

\title[A new approach to inverse Sturm-Liouville problems]
{A new approach to inverse Sturm-Liouville problems based on point interaction}

\author[Min Zhao]{Min Zhao}
\address{School of Mathematics and Statistics, Shandong University, Weihai 264209, P.\ R.\ China}
\email{zhaomin215@mail.sdu.edu.cn}\

\author[Jiangang Qi]{Jiangang Qi}
\address{School of Mathematics and Statistics, Shandong University, Weihai 264209, P.\ R.\ China}
\email{qjg816@163.com}\

\author[Xiao Chen]{Xiao Chen$^\dag$}
\address{School of Mathematics and Statistics, Shandong University, Weihai 264209, P.\ R.\ China}
\email{chenxiao@sdu.edu.cn}\


\begin{abstract}
In the present paper, motivated by point interaction, we propose a new and explicit approach to inverse Sturm-Liouville eigenvalue problems under Dirichlet boundary.
More precisely, when a given Sturm-Liouville eigenvalue problem with the unknown integrable potential interacts with $\delta$-function potentials, we obtain a family of perturbation problems, called \emph{point interaction models} in quantum mechanics.
Then, only depending on the first eigenvalues of these perturbed problems, we define and study \emph{the first eigenvalue function}, by which the desired potential can be expressed explicitly and uniquely. As by-products, using the analytic function theoretic tools, we also generalize several fundamental theorems of classical Sturm-Liouville problems to measure differential equations with non-trivial positive integrable weights.

\vspace{04pt}


\end{abstract}

\begin{thanks}{$\dag$ the corresponding author, chenxiao@sdu.edu.cn.}
\end{thanks}

\maketitle

\section{Introduction and problem statement}\label{sec:intro}

It is well known that, in quantum mechanics, the eigenvalues (i.e., point spectrum) are used to describe the energies of some particles under certain state, and also are the only observables.
Hence, it is especially significant to study how to determine and reconstruct a quantum system by eigenvalues.
Since Sturm-Liouville problems (abbreviated as SLPs or S-L problems), which characterize the stationary quantum systems, usually depend on their potentials and weights, it is fairly easy to understand that the classical inverse problems for SLPs mainly focus on the studies of how to determine the potential uniquely using suitable spectral data such as eigenvalues.
So the inverse SLPs have always been the subjects of intense scholarly research, and have numerous applications to many fields in mathematics and natural science (cf. \cite{ GGKM1967, LY2006, FY2008, FWW2016}, etc).

The first result in the inverse SLPs, which is due to V.\ A.\ Ambarzumian in 1929, showed that the zero potential for a special class of SLPs can be recovered by one spectra.
In 1946, G.\ Borg further proved the famous foundational theorem that the potentials of SLPs can be determined by two spectra.
Afterwards N.\ Levinson in 1949 further improved Borg's result as well as simplified the proof. Thereby, the cornerstone of the inverse Sturm-Liouville theory had been laid down.
From then on, a large number of literatures have been devoted to this direction, wherein a variety of methods of dealing with these problems have come into being, for example, the classical Gelfand-Levitan theory, the Simon theory, the approach proposed by Remling, the boundary control method, the Hochstadt-Liebermanthe method, the geometric method, the Riesz-basis method, the spectrum-like function and so on. For more details about these methods above, the reader may refer to a review article \cite{AM2010} that briefly introduces the first four methods above,  and also may see \cite{FY2008, FWW2016, PT1987, KK1997, Gesz2007, YB2020, CQ2021} which contains the other methods.

Among the existing methods, for recovering the unknown potentials on the whole interval, we need at least two complete sets of eigenvalues that are pre-given.
Unfortunately, in practice, we generally are only able to obtain incomplete spectral data, such as especially the finite eigenvalues.
So more and more people have to try to use  the partial spectral data to recover the potentials (cf. \cite{WX2009, KNT2017, Krav2022}, etc).
It is worth noticing that, in these cases, we often cannot make sure that a potential can be recovered uniquely.
Thereupon,  we naturally begin to consider inverse problems of characterizing the set consisting of all eligible potentials reconstructed by given incomplete spectral data, and the optimization problems of calculating the infimum of the $L^1$-norm of all such potentials as well as finding the optimal elements attaining the infimum (cf. \cite{LY2006, CQ2021, QC2016, WZZ2020, GQ2020}, etc).
Here we will be devoted to recover the potentials uniquely from the first eigenvalues.

The main idea of our method comes from the {\bf $\delta$-point interaction} in quantum physics.
The $\delta$-point interaction model, which is a system of particles of equal masses in one dimension interacting through a $\delta$-function potential such as $-\frac{\hbar^2}{2m}\frac{\rd^2}{\rd x^2}+g\delta(x)$, is one of the exactly solvable many-body models in quantum mechanics.
The $\delta$-function potential is a special case of a large family of point interactions in one dimension.
Historically, the first influential paper \cite{KP1931} on the $\delta$-point interaction model was that by Kronig and Penney in 1931.
The Kronig-Penney model has become a standard reference model in solid state physics.
To this day, a lot of physicists and mathematicians, including the Nobel Laureate C.\ N.\ Yang \cite{Yang1967}, have already made a substantial contribution to this field. For more details about point interaction, the reader may refer to \cite{AGHH1988, BMK2023} and the references therein.

In this paper,  we mainly consider the S-L Dirichlet eigenvalue problem $({\bf E}_q)$ defined as:
\begin{equation}\label{eqn:main-probl}
-y''(x)+q(x)y(x)=\lambda w(x)y(x) \text{ on } [0,1],\ y(0)=0=y(1),
\end{equation}
where $q,\ w\in L^1[0,1]$ with $w>0$ a.e. on $(0,1)$, and $\lambda$ is the spectral parameter.
Denote by $\lambda_1(q)$ the first eigenvalue of \eqref{eqn:main-probl}.
The potential $q$ is unknown, and note that for any fixed $\lambda\in\mathbb{R}$, the set
$S(\lambda):=\{q :\ \lambda_1(q)=\lambda,\ q\in L^1[0,1]\}$ is  an infinite-dimensional submanifold of $L^1[0,1]$, see \cite[p.68]{PT1987}.
Therefore, it is impossible to recover $q$ uniquely by $\lambda_1(q)$.

At first, for any given $t\in(0,1)$ and $r\in [0,\epsilon_0]$ with a small enough number $\epsilon_0>0$, we let the original system \eqref{eqn:main-probl} interact with a $\delta$-function potential in the following manner:
\begin{equation}\label{eqn:pt-interact-probl}
-y''(x)+[q(x)-r\delta(x-t)]y(x)=\lambda w(x)y(x), \ y=y(x)\text{ on }(0,1),\ y(0)=0=y(1),
\end{equation}
where $q,\ w$ are just the ones in \eqref{eqn:main-probl}, the positive number $r$ is called {\bf coupling constant} that stands for the intensity of the interaction, and $\delta(x-t)$ is the {\bf Dirac $\delta$-function} at $t\in(0,1)$ defined by 
\begin{equation}\label{eqn:delta-funct}
\delta(x-t)=
\begin{cases}
\infty, &x=t,\\
          0, &x\neq t
\end{cases}
\ \text{and}\ \int_I\delta(x-t)\, dx=1,\ \ \forall\, I\subset[0,1]
\ \text{and}\ t\in I.
\end{equation}
Note that $$\int_0^1 f(x)\delta(x-t)\,{\rm d}x=\int_{t-\epsilon}^{t+\epsilon} f(x)\delta(x-t)\,{\rm d}x=f(t),$$ for any continuous function $f$ and $0<\epsilon<\min\{t,1-t\}$.
As  $t$ (resp. $r$) runs over all real numbers in $(0,1)$ (resp. $[0, \epsilon_0)$),  we obtain a family of {\bf perturbation problems} $({\bf E}^{t,r}_q)$, defined as \eqref{eqn:pt-interact-probl}, of the original problem $({\bf E}_q)$. This kind of equations also can be found in \cite{Mana2016, BMK2023}.

In general, the coupling constant $r$ can be replaced by a continuous function $r(t)$ with respect to the interaction position $t$, and we may call $r(t)$ {\bf coupling coefficient function}.
This sort of perturbation problem also can be viewed as a {\bf moving point interaction} model \cite{Pos2007}.
Notice that the equations and their solutions in both \eqref{eqn:main-probl} and \eqref{eqn:pt-interact-probl} can be explained in the senses of both ordinary differential equation with distributions \cite{Mana2016, BMK2023, CQ2021} and measure differential equation \cite{MZ2013, ZWMQX2018}.

Next, suppose that the first eigenvalues $\lambda(t,r;q)$ of the perturbation problem $({\bf E}^{t,r}_q)$ are known for any $(t, r)\in (0,1)\times [0,\epsilon_0]$, which means that the first energy eigenvalues of all perturbation systems can be observed. And then it will be proved that $\lambda(t,r;q)$, which here is named the {\bf first eigenvalue function} of $({\bf E}_q)$, is continuous and differentiable with respect to $(t, r)\in(0,1)\times[0,\epsilon_0]$.

Lastly,  we uniquely and directly reconstruct the desired potential $q$ of the main problem $({\bf E}_q)$ by
\begin{equation}\label{eqn: q-recover}
q(x)=\frac{\varphi''_0(x)}{\varphi_0(x)}+\lambda_1 w(x),
\end{equation}
where $\lambda_1=\lambda(x,0;q)$ and $\varphi_0(x)=\sqrt{-\frac{\partial \lambda(x,0;q)}{\partial r}}$ on $(0,1)$. Note that $\lambda(x,0;q)$ always equals to $\lambda_1(q)$ since the problem $({\bf E}^{t,0}_q)$ is actually $({\bf E}_q)$. See Theorem~\ref{thm:main-thm}.

In \cite{CQ2021}, the inverse Dirichlet problem with single Dirac weight also can be regarded as a homogeneous S-L equation with $\delta$-point interaction, if the Dirac weight, shifted to the left hand side of the equation, is viewed as a Dirac potential, and then the complete spectral data $\lambda(t)$ is exactly coupling coefficient function mentioned above. Consequently, the inverse eigenvalue problem in \cite{CQ2021} can be transformed to an inverse non-eigenvalue problem with point interaction.

This paper is organized as follows.
In Section~\ref{sec:prelim}, as a preliminary, we will introduce some basic terminologies, notations and useful facts about measure differential equations, and also prove weak$^*$ continuity of eigenvalues in measures for SLPs with measure potentials and positive integrable weights (i.e., Theorem~\ref{thm:cont-eigenv}).
In Section~\ref{sec:1st-eigen-funct}, we will define and study the first eigenvalue function.
Section~\ref{sec:potential-recovery} is the core part,  wherein we prove the uniqueness (but also reconstruction) result of recovering potential by the first eigenvalue function (i.e., Theorem~\ref{thm:main-thm}).
In Appendix~\ref{sec:app}, applying the holomorphic technique, we generalize several fundamental theorems of classical Sturm-Liouville problems to measure differential equations with measure potentials and positive integrable weights, for example, weak$^*$-continuity of eigenvalues in measure (i.e., Theorem~\ref{thm:cont-eigenv}), the geometric multiplicity of eigenvalues (i.e., Theorem~\ref{thm:eigenv-simply}) and the number of zeros of eigenfunctions (i.e., Theorem~\ref{thm:zero-eigenf}).  Notice that the proofs are quite different from those \cite{MZ2013} in the case that the weight is the constant $1$.

\bigskip

\section{Notations and preliminary}\label{sec:prelim}
\medskip

In this section, we recall some materials on measure theory and measure differential equation, and also prove several related results we need later (cf. \cite{Fol1999, MZ2013, ZWMQX2018} and the references therein).

Within this paper, we denote by $\BR$ (resp., $\BN$ and $\BZ^+$)  the set of all real numbers (resp., non-negative integers and positive integers).
Let $\BK$ be $\BR$ or $\BC$.
Denote by $C_\BK[0,1]$ (resp. $L^1_\BK[0,1]$) the Banach space generated by all $\BK$-value continuous (resp. Lebesgue integrable) functions on $[0,1]$ under the supremum norm $\norm{\cdot}_\infty$ (resp. $L^1$-norm $\norm{\cdot}_{L^1}$). For any subset $U$ of $[0,1]$, denote by $\chi_U$  the characteristic function on $U$, defined by 
$$\chi_U(x)=1,\ \forall x\in U;\quad \chi_U(x)=0,\ \forall x\in [0,1]\setminus U.$$
Let $AC_\BK[0,1]$ be the space consisting of all {\bf absolutely continuous (a.c.)} $\BK$-value functions on $[0,1]$.

For a real function $\mu: [0,1]\rightarrow \mathbb{K}$,  the {\bf total variation} of $\mu$ on $[0,1]$ is defined as
$$\norm\mu_{\bf V}=\sup \left\{\sum\limits_{i=0}\limits^{n-1}|\mu(x_{i+1})-\mu(x_{i})| :\ 0=x_{0}<\cdots<x_{n}=1,n\in\mathbb{N}\right\}.$$
The space of normalized $\BK$-value functions of {\bf bounded variation (NBV)} on $[0,1]$ is defined as
$$\CM^\BK_0:=\left\{\mu: [0,1]\rightarrow \BK :\ \mu(0^+)=0,\ \mu(x^+)=\mu(x),\ \forall x\in (0,1),\  \norm{\mu}_{\bf V}<\infty\right\}.$$
Here $\mu(x^+):=\lim_{s\rightarrow x^+}\mu(s)$ is the right-limit for any $x\in [0,1)$, and the word ``normalized" means the condition $\mu(0+)=0$. 
Without this condition, $\mu$ is just called a function of {\bf bounded variation (BV)}.
Every $\mu\in\CM^\BK_0$ induces a unique $\BK$-value Borel measure $\tilde\mu$ on $[0,1]$ such that $\mu(x)=\tilde\mu((0,x])$ for any $x\in(0,1]$, and vice versa. 
Namely, this $\mu$ is the distribution function of the measure $\tilde\mu$. Hence, throughout this paper, we identify $\mu$ with $\tilde\mu$ if without confusion, and also view $\CM^\BK_0$ as the Banach space over $\BK$ consisting of all $\BK$-value Borel measures on $[0,1]$ under the norm $\norm{\mu}_{\bf V}$.

Obviously, the Lebesgue measure on $[0,1]$ corresponds to $\ell(s)=s\in\CM^\BR_0$, i.e., ${\rm d}\ell={\rm d}s$, and every $f\in L^1_\BK[0,1]$ defines an a.c. measure  in $\CM^\BK_0$ by the Riemann-Stieltjes integral
\begin{equation}\label{eqn: ac-mes}
\mu_{f}(x):=\int_{[0,x]}f(s)\, {\rm d}s,\quad \text{on }[0,1], \text{ i.e., } \mu_f\in AC_\BK[0,1] \text{ and } d\mu_f=fds.
\end{equation} 
In particular, the {\bf Dirac measure} $\delta_t$ (also called the point mass, or Heaviside function) at $t\in(0,1)$, defined as
\begin{equation}\label{eqn:delta-mes}
\delta_t(x)=\left\{
\begin{array}{cl}
0, &  x\in[0,t), \\
1,  &  x\in[t,1], \\
\end{array} \right.
\end{equation}
is a discrete measure in $\CM^\BK_0$. The Dirac $\delta$-function  $\delta(x-t)$ defined by \eqref{eqn:delta-funct} is exactly the Radon-Nikodym derivative of the Dirac measure $\delta_t$, i.e., ${\rm d}\delta_t=\delta(x-t){\rm d}x$ customarily.

Moreover, the famous Riesz representation theorem tells us that $(\CM^\BK_0, \norm{\cdot}_{\bf V})$ is isomorphic to the dual space $(C_\BK[0,1], \norm{\cdot}_\infty)^*$. 
In fact, for every $\mu\in\CM^\BK_0$, we can define a bounded linear functional on $C_\BK[0,1]$ by the Riemann-Stieltjes integral
$$\mu^*(u):=\int_{[0,1]} u(x)\,{\rm d}\mu(x),\quad \forall u\in C_\BK[0,1],$$
which induces the weak$^*$ ($w^*$-) topology on $\CM^\BK_0$.

\begin{defn}\cite[Definition~2.2]{MZ2013}\label{defn:w-star-top}
For any  $\mu_{n},\ \mu_0\in \CM^\BK_{0}$ and $n\in \mathbb{N}$, we say that  $\mu_{n}$ is weakly$^*$ convergent to $\mu_0$, denoted as $\mu_{n} \overset{w^*}{\longrightarrow}\mu_0$ in $(\CM^\BK_{0},w^{\ast})$, if
$$\lim\limits_{n\rightarrow\infty}\int_{[0,1]}u(t)\,{\rm d}\mu_{n}(t)=\int_{[0,1]}u(t)\,{\rm d}\mu_0(t),\quad \forall u\in C_\BK[0,1].$$
\end{defn}

\begin{rem}\label{rem:L1-approx}
$(i)$ It is well known that $\CM^\BK_0$ can be decomposed into the direct sum $\CM^\BK_{ac}\oplus\CM^\BK_{sc}\oplus\CM^\BK_d$. 
Here $\CM^\BK_{ac}$ is the a.c. part $AC_\BK[0,1]$. 
The second component $\CM^\BK_{sc}$ is the {\bf singular continuous (s.c.)} part, the one in which is continuous and has zero Radon-Nikodym derivative almost everywhere on $[0,1]$.   
The third one $\CM^\BK_d$ is the discrete (also called pure discontinuous or complete singular) part that is 
$$\left\{\sum_{a\in\CA} m_a\delta_a:\  \CA\text{ is at most countable subset of }[0,1],\ m_a\in\BK, \text{ and }\sum_{a\in A} \abs{m_a}<+\infty\right\},$$
where every element is a step function having at most countable discontinuous points.
All these subspaces of $\CM^\BK_0$ are norm-closed. For more details, refer to \cite[Page~106]{Fol1999},  \cite[Lemma~2.2]{ZWMQX2018} and the references therein.

\smallskip
\noindent
$(ii)$ Note that for any $ \mu \in \CM^\BK_d$, there exists $\{f_n\}_{n\in\BZ^+}$ in $L^1_\BK[0,1]$ such that $\mu_{f_n}\overset{w^*}{\longrightarrow}\mu$.
Indeed, for any $\delta_t\in\CM_d$, set $\varepsilon:=\min\{t,1-t\}$, it is apparent that $\frac{n}{2\varepsilon}\chi_{[t-\frac{\varepsilon}{n},t+\frac{\varepsilon}{n}]}$ is weakly$^*$ convergent to $\delta_t$ in $\CM^\BK_0$. 

\smallskip
\noindent
$(iii)$Note that $\CM^\BR_{ac}$ is dense in $(\CM^\BR_0,w^*)$, that is, for any $ \mu \in \CM^\BR_0$, there exists $\{f_n\}_{n\in\BZ^+}$ in $L^1_\BR[0,1]$ such that $\mu_{f_n}\overset{w^*}{\longrightarrow}\mu$.
Indeed, since every NBV function can be pointwisely approximated by a sequence of a.c. functions, it follows from  \cite[proposition~7.19]{Fol1999} as well as the continuity of the s.c. functions, that $\CM^\BR_{sc}$ is contained in the weak$^*$ closure of $\CM^\BR_{ac}$. Therefore, by the statements $(1)$ and $(2)$ above, we obtain the weak$^*$ denseness of $\CM^\BR_{ac}$ in $\CM^\BR_0$.
\end{rem}

Consider the initial value problem of the {\bf measure differential equation (MDE)} defined as
\begin{equation}\label{eqn:mde}
\left\{\aligned
&-{\rm d}y^{\bullet}(x)+y(x)\,{\rm d}\mu(x)=0,\ x\in [0,1],\\
&(y(0),y^{\bullet}(0))=(y_0, z_0)\in\BK^2,
\endaligned\right.
\end{equation}
where $\mu\in \CM^\BK_{0}$. The solution $y(x)\in C[0,1]$ and its generalized derivative $y^{\bullet}(x)$ are defined by the following integral system:
\begin{equation}\label{eqn:solution-mde}
\left\{\aligned
&y(x)=y_0+\int_{[0,x]} y^\bullet(t)\,{\rm d}t,\ x\in [0,1],\\
&y^\bullet(x)=\left\{\aligned
& z_0,\ x=0,\\
& z_0+\int_{[0,x]} y(t)\,{\rm d}\mu(t),\ x\in(0,1].
\endaligned\right.
\endaligned\right.
\end{equation}
Here the first integral is Lebesgue integral, and the second one is Riemann-Stieltjes integral.  
For any $(y_0, z_0)\in\BK^2$, there exists a unique solution $y$ for the initial value problem \eqref{eqn:mde}. The derivative $y^\bullet$ is a real BV function on $[0,1]$, while $y$ is actually absolutely continuous. 
Moreover, the derivative $y^\bullet$ is just the classical right-derivative $y'_+(x)$ on $(0,1)$, and $y^\bullet=y'$ a.e. on $[0,1]$. 
For more details, refer to \cite[Definition~3.1--Corollary~3.4]{MZ2013}, \cite[Eq.\,(1.3)-(1.4)]{ZWMQX2018} and the references therein.

The problem $({\bf E}^{t,r}_q)$ can be regarded as a Dirichlet eigenvalue problem of a MDE \eqref{eqn:mde} with $\mu=\mu_q-r\delta_t-\lambda\mu_w\in\CM^\BR_0$ as well as Dirichlet boundary $y(0)=y(1)=0$, and so its solution can be explained as in \eqref{eqn:solution-mde}. On the other hand, since the problem $({\bf E}^{t,r}_q)$ is actually a Sturm-Liouville eigenvalue problem with distribution coefficients, similarly to \cite{CQ2021}, its solution $y$ also can be determined by 
$$y\in \left\{z:\ z\in AC_\BR[0,1],\ z'\in [0,t)\cup(t,1],\ \exists z'(t^\pm)\right\}$$
and
\begin{equation}\label{eqn:solution-SL-distrib}
\left\{\aligned
&-y''(x)+q(x)y(x)=\lambda w(x)y(x),\ t\neq x\in [0,1],\\
&y'(t^-)-y'(t^+)=ry(t)\\
&y(0)=0=y(1).\\
\endaligned\right.
\end{equation}

As the end, we consider the Dirichlet eigenvalue problem $({\bf E}_\mu)$ of MDE as follows:

\begin{equation}\label{eqn:mu-probl-mde}
\left\{\aligned
&-{\rm d}y^{\bullet}(x)+y(x)\,{\rm d}\mu(x)=\lambda w(x)y(x)\,{\rm d}x,\ x\in [0,1],\\
&y(0)=0=y(1),
\endaligned\right.
\end{equation}
where $\mu\in \CM^\BR_{0}$, $w\in L^1[0,1]$, $w>0$ a.e. on $[0,1]$, and $\lambda$ is the spectral parameter.
Clearly, the original problem $({\bf E}_q)$ and perturbation problem $({\bf E}^{t,r}_q)$ are exactly the problems $({\bf E}_{\mu_q})$ and $({\bf E}_{\mu_q-r\delta_t})$, respectively.

The authors in \cite[Theorem~1.3]{MZ2013} used Pr\"ufer transformation to prove the weak$^*$ continuity of eigenvalues in measures for the problem $({\bf E}_\mu)$ with $w\equiv 1$. 
Normally, as a natural generalization of \cite[Theorem~1.3]{MZ2013}, Theorem~\ref{thm:cont-eigenv} also should hold. 
However, we do not find it explicitly proved in the literature, so we will give the complete argument in Appendix~\ref{sec:app} for benefit of the reader. 
It is worth mentioning that, instead of Pr\"ufer transformation that seems not to work well for Theorem~\ref{thm:cont-eigenv}, we will apply the analytic function theory to present a new proof. 
Moreover, we in Appendix~\ref{sec:app} also generalize some classical theorems about eigenvalue and eigenfunction to the Dirichlet eigenvalue problem $({\bf E}_\mu)$ of MDE with positive integrable weights.

\begin{thm}\label{thm:cont-eigenv}
Let $\lambda_m(\mu)$ denote the $m$-th eigenvalue of the problem $({\bf E}_\mu)$ defined as in \eqref{eqn:mu-probl-mde} with the measure $\mu\in\CM^\BR_{0}$. Then $\lambda_m(\mu)$ is continuous in the measure $\mu\in(\mathcal{M}_{0},w^{\ast})$.
\end{thm}

\bigskip
\section{The first eigenvalue function $\lambda(t,r)$}\label{sec:1st-eigen-funct}
\medskip

In this section, we  propose and study a new concept, called \emph{the first eigenvalue function of a SLP}, which is crucial for our paper.  

Recall that the main problem $({\bf E}_q)$ can be written as the Dirichlet problem $({\bf E}_{\mu_q})$ of a MDE of the form
\begin{equation}\label{eqn:main-probl-mde}
-dy^\bullet(x)+y(x)q(x)\, {\rm d}x=\lambda w(x)y(x)\, {\rm d}x,\ y(0)=0=y(1),\ x\in[0,1].
\end{equation}
where $q$ (i.e., the a.c. measure $\mu_q$ and ${\rm d}\mu_q=q(x)\,{\rm d}x$) is unknown.

In this paper, to recover the unknown $q$, we try to adopt a new method, that is, we use the first eigenvalues of a family of perturbation problems to determine the desired potential. More precisely, we consider the perturbed S-L problems $({\bf E}^\mu_q)$ written as
\begin{equation}\label{eqn:perturb-probl}
-dy^\bullet(x)+y(x)(q(x)\, {\rm d}x+d\mu(x,t))=\lambda w(x)y(x)\, {\rm d}x,\ y(0)=0=y(1),\ x\in[0,1].
\end{equation}
where $\mu(x,t)\in\CM^\BR_0$ is the perturbation term with a parameter $t\in[0,1]$. 

Denote by $\lambda(\mu;q)$ the first eigenvalue of $({\bf E}^\mu_q)$. 
Then, it is quite reasonable that we ask if the $q$ can be uniquely determined by $\lambda(\mu;q)$. Apparently, it is not always successful for all perturbations $\mu\in\CM^\BR_0$.

For example, if $d\mu_0=ctw(x)\,{\rm d}x$ is chosen to be the perturbation in \eqref{eqn:perturb-probl}, where $c\in\mathbb{R}$ is a constant and $t\in[0,1]$, then
\begin{equation}\label{eqn:sick-perturb}
\lambda(\mu_0;q)=\lambda_1(q)+ct,
\end{equation}
where $\lambda_1(q)$ is the first eigenvalue of the problem $({\bf E}_{\mu_q})$ (i.e., the main problem $({\bf E}_q)$).
Unfortunately, by Eq.\,\eqref{eqn:sick-perturb}, we can see that, for any two distinct potentials $q_1,q_2\in S(\lambda_0)$ mentioned in Section~\ref{sec:intro}, it follows that $\lambda(\mu_0;q_1)=\lambda(\mu_0;q_2)=\lambda_0+ct$. 
In other words,  we can not distinguish different potentials through the first eigenvalues $\lambda(\mu_0)$ under such perturbation $\mu_0$ above. 
For such reason, the perturbation $\mu$ in the problem $({\bf E}^\mu_q)$  should be chosen skillfully. 

Inspired by $\delta$-point interaction in quantum mechanics, we take the point interaction $\mu=-r\delta_t$ as the perturbation term in \eqref{eqn:perturb-probl}, that is, we are going to consider the perturbed SLPs $({\bf E}_q^{t,r})$.
For a given problem $({\bf E}_q^{t,r})$ with $(r,t)\in (0,1)\times [0,\epsilon_0]$,  we assume that the potential $q\in L^1_\BR[0,1]$ is unknown, and meanwhile the first eigenvalue is available. 
By moving the position $t$ of interaction along the open interval $(0,1)$ as well as adjusting the intensity $r$, we gain a family of the first eigenvalues, denoted by $\lambda(t,r;q)$, of $({\bf E}_q^{t,r})$ for all $t\in(0,1)$ and $r\in [0,\epsilon_0]$. 
In the remainder of this section, we will prove the continuity and differentiability of $\lambda(t,r;q)$ in $(t,r)\in (0,1)\times [0,\epsilon_0]$, and so view $\lambda(t,r;q)$ as a function related to $({\bf E}_q)$.

\begin{defn}\label{defn:1st-eigenv-funct}
We call the $\lambda(t,r;q)$ above the {\bf first eigenvalue function} of the problem $({\bf E}_q)$ defined by \eqref{eqn:main-probl}.
\end{defn}

\noindent
{\bf Notation: }for a given problem $({\bf E}_q)$,  the potential $q$ has already been fixed, although it is unknown. Hence, for convenience, thereinafter {\bf we sometimes simply denote by $\lambda(t,r)$ the first eigenvalue function when there is no risk of confusion}.

\begin{thm}\label{thm:cont-1st-eigenv-funct}
Let $\lambda(t,r)$ be the first eigenvalue function of $({\bf E}_q)$.
Then $\lambda(t,r)$ is continuous on $(t,r)\in [0,1]\times [0,+\infty)$,
and the corresponding eigenfunction  $\Phi(x,\lambda(t,r))$ does not change its sign on $(0,1)$. Furthermore, it holds that
$$
\lambda(t,r)<\lambda(0,r)=\lambda(1,r)=\lambda(t,0)=\lambda_1,\quad (t,r)\in(0,1)\times(0,+\infty),
$$
where $\lambda_1$ is the first eigenvalue of $({\bf E}_q)$.
\end{thm}
\begin{prf}. Clearly, $\lambda(t,0)=\lambda_1$ for all $t\in[0,1]$.
If $t=0$(reps. $t=1$), the perturbation term $-r\delta(x-0)$ (reps. $-r\delta(x-1)$)
has no influence to the original problem  $({\bf E}_q)$, and hence $\lambda(0,r)=\lambda(1,r)=\lambda_1$.

Since $$\int^1_0[r_n\delta(x-t_n)-r\delta(x-t)]f(x)\,{\rm d}x=r_nf(t_n)-rf(t)\to 0,\ n\to+\infty,\ \forall f\in C_\BR[0,1],$$
we have that $\lambda(t,r)$ is continuous on $(t,r)$ by Theorem \ref{thm:cont-eigenv}.

For fixed $(t,r)\in(0,1)\times(0,+\infty)$, from Remark \ref{rem:L1-approx}$(ii)$, choose $p_n(x)\in L^1[0,1]$ such that
$$
p_n(x)\geqslant 0,\ \int^1_0 p_n(x)\,{\rm d}x=1,\  \mu_{p_n} \stackrel{w^{\ast}}{\longrightarrow} \delta_t,\ n\to+\infty.
$$
Then $\mu_{q-rp_n}\stackrel{w^{\ast}}{\longrightarrow} \mu_q-r\delta_t$ as $n\to+\infty$, and hence $\lambda_1(p_n)\to \lambda(t,r)$ as $n\to+\infty$ by Theorem \ref{thm:cont-eigenv}, where $\lambda_1(p_n)$ is the first eigenvalue of the problem
$$({\bf E}_{\mu_{q-rp_n}}):\  -y''+[q(x)-rp_n(x)]y=\lambda w y,\ y(0)=0=y(1). \ \text{ See \eqref{eqn:mu-probl-mde}.}$$
According to the monotonicity of eigenvalues with respect to potential (\cite[Theorem4.9.1]{Zettl2005}), the inequality $q(x)-rp_n(x)\leqslant q(x)$ implies $\lambda_1(p_n)\leqslant\lambda_1$, and hence $\lambda(t,r)\leqslant\lambda_1$.

Now, we prove that $\lambda(t,r)<\lambda_1$ if $(t,r)\in(0,1)\times(0,+\infty)$. Let $\varphi(x,\lambda)$ and $\psi(x,\lambda)$ be solutions
of $-y''+qy=\lambda w y$ such that
\begin{equation}\label{eqn:phi-psi-ini-value}
 \varphi(0)=0, \ \varphi'(0)>0;\  \psi(1)=0, \ \psi'(1)<0.
\end{equation}
Denoted $\Delta(\lambda)$ by the  Wronskian $W[\varphi,\psi](x)$ of $\varphi$ and $\psi$, i.e.
\begin{equation}\label{eqn:wronski}
W[\varphi,\psi](x)=\varphi(x,\lambda)\psi'(x,\lambda)-\varphi'(x,\lambda)\psi(x,\lambda)\equiv \Delta(\lambda).
\end{equation}
We remark that
the choice of $\varphi$ (resp. $\psi$) is not unique, but this doesn't matter
in the following discussion. Let  $\Phi(x,\lambda(t,r))$ be the first eigenfunction of $({\bf E}^{t,r}_q)$, that is
\begin{equation}\label{eqn:Phi-eqn}
-\Phi''+[q-r\delta(x-t)]\Phi=\lambda(t,r) w \Phi,\ \Phi(0)=0=\Phi(1),
\end{equation}
or equivalently, by \eqref{eqn:solution-SL-distrib},
\begin{equation}\label{eqn:defn-Phi}
\left\{\aligned
&-\Phi''+q \Phi=\lambda(t,r) w \Phi,\ t\not=x\in(0,1),\\
& \Phi'(t-0,\lambda(t,r))-\Phi'(t+0,\lambda(t,r))=r\Phi(t,\lambda(t,r)),\\
& \Phi(0,\lambda(t,r))=0=\Phi(1,\lambda(t,r)),
\endaligned\right.
\end{equation}
where $\Phi'=\frac{\partial\Phi(x,\lambda(t,r))}{\partial x}$.
Then $\Phi$ can be chosen as
$$
\Phi(x,\lambda(t,r))=\left\{\aligned
&\varphi(x,\lambda(t,r)),\ &x\in[0,t],\\
&c(t,r)\psi(x,\lambda(t,r)),\  &x\in[t,1],
\endaligned\right.
$$
where $c(t,r)=\frac{\varphi(t,\lambda(t,r))}{\psi(t,\lambda(t,r))}$. Since $\lambda(t,r)\leqslant\lambda(0)=\lambda_1$ and
the first eigenfunction of $({\bf E}_q)$ associated to $\lambda_1$ is positive on $(0,1)$ (\cite[Theorem 4.3.1(6)]{Zettl2005}),
it follows from the oscillation theory of linearly second order differential equations (\cite[Theorem 2.6.2]{Zettl2005})
that
$\varphi(x,\lambda(t,r))>0$, $\psi(x,\lambda(t,r))>0$ for $x\in(0,1)$,  and hence $c(t,r)>0$, which means that
$\Phi(x,\lambda(t,r))>0$ on $(0,1)$ for all $(t,r)\in[0,1]\times[0,+\infty)$.

Let $\varphi_0$ be the corresponding eigenfunction of $({\bf E}_q)$ associated to $\lambda_1$, that is
\begin{equation}\label{eqn:phi_0-Eq}
-\varphi''_0+q\varphi_0=\lambda_1 w\varphi_0, \ \varphi_0(0)=0=\varphi_0(1).
\end{equation}
Multiplying two sides of the equations \eqref{eqn:Phi-eqn} and \eqref{eqn:phi_0-Eq} by $\varphi_0$ and $\Phi$ respectively, and then
integrating their difference over $[0,1]$, we find that
$$
r\varphi_0(t)\Phi(t,\lambda(t,r))=-(\lambda(t,r)-\lambda_1)\int^1_0w(x)\varphi_0(x)\Phi(x,\lambda(t,r))\,{\rm d}x,
$$
which, together with the positivity of  both $\varphi_0$ and $\Phi(x,\lambda(t,r))$, implies that $\lambda(t,r)<\lambda_1$
for $(t,r)\in(0,1)\times(0,+\infty)$.
\end{prf}

The following theorem  gives the necessary  and sufficient condition for a function $\lambda(t,r)$ being the first eigenvalue function of $({\bf E}_q)$.

\begin{thm}\label{thm: equiv-1st-eigenv}
Let $\varphi(x,\lambda)$ and $\psi(x,\lambda)$ be defined as in \eqref{eqn:phi-psi-ini-value}.
Then $\lambda(t,r)$ is the first eigenvalue function of $({\bf E}_q)$ if and only if
\begin{equation}\label{eqn:equiv-1st-eigenv}
 r\varphi(t,\lambda(t,r))\psi(t,\lambda(t,r))=-\Delta(\lambda(t,r)),\ r>0,\ t\in(0,1),
\end{equation}
where $\Delta(\lambda)=W[\varphi,\psi]$ is the Wronskian \eqref{eqn:wronski} of $\varphi, \psi$.
\end{thm}

\noindent
\begin{prf} 
Note that  the eigenfunction  $\Phi(x,\lambda(t,r))$ of $({\bf E}^{t,r}_q)$ can be chosen as
\begin{equation}\label{eqn:eigenf-perturb-probl}
\Phi(x,\lambda)=\left\{\aligned
&\varphi(x,\lambda),\ &x\in[0,t]\\
&c(t,r)\psi(x,\lambda),\  &x\in[t,1],
\endaligned\right.\quad
c(t,r)=\frac{\varphi(t,\lambda)}{\psi(t,\lambda)},
\end{equation}
where $\lambda$ simply denotes $\lambda(t,r)$. It follows from the second equation in \eqref{eqn:defn-Phi} that
$$
\varphi'(t,\lambda)-c(t,r)\psi'(t,\lambda)=r\varphi(t,\lambda).
$$
As a result, it holds that
$$
r\varphi(t,\lambda)=\varphi'(t,\lambda)-\frac{\varphi(t,\lambda)}{\psi(t,\lambda)}\psi'(t,\lambda),
$$
which implies Eq.\,\eqref{eqn:equiv-1st-eigenv}. Since every step in the above discussion can be coversed,
we can know that $\lambda(t,r)$ satisfying Eq.\,\eqref{eqn:equiv-1st-eigenv}
implies $\lambda(t,r)$ being the first eigenvalue function of $({\bf E}_q)$.
\end{prf}

\begin{rem}
Note that $\lambda(t,r)<\lambda_1$ for $(t,r)\in(0,1)\times (0,+\infty)$ and
$\varphi(t,\lambda)>0$ on $t\in(0,1]$ for $\lambda<\lambda_1$.
If we replace $\psi$ by $\widehat{\psi }$, where
$$
\widehat{\psi}(x)=\varphi(x)\int^1_x\frac{ds}{\varphi^2(s,\lambda)},
$$
then Eq.\,\eqref{eqn:equiv-1st-eigenv} can be written as
\begin{equation}\label{eqn:int-eqn}
 r\varphi^2(t,\lambda(t,r))\int^1_t\frac{ds}{\varphi^2(s,\lambda(t,r))}=1,\ r>0,\ t\in(0,1)
\end{equation}
if $W[\varphi,\widehat{\psi}](x)\equiv-1$. We should point out that for unknown potential
$q$, if we know the corresponding first eigenvalue function $\lambda(t,r)$, the equation \eqref{eqn:int-eqn}
is a kind of integral equation involving the unknown function $\varphi$. Clearly, if the unknown  $\varphi$
is able to be solved from \eqref{eqn:int-eqn}, then $q$ can be expressed by
$$
q(x)=\frac{\varphi''(x,\lambda(t,r))}{\varphi(x,\lambda(t,r))}+\lambda(t,r)w(x).
$$
But it is very difficult to solve $\varphi$ directly due to its nonlinearity. For this reason, we will
seek for other efficient methods to recover the potential $q$.
\end{rem}

To this end, we further obtain the differentiability of the first eigenvalue function $\lambda(t,r)$ and the corresponding normalized eigenfunction. 

\begin{thm}\label{thm:diff-1st-eigenv-funct}
If $\lambda(t,r)$ is  the first eigenvalue function of $({\bf E}_q)$, then one has

\noindent
$(i)$ $\lambda(t,r)$ is continuous differentiable with respect to $(t,r)\in(0,1)\times[0,+\infty)$ and
\begin{equation}\label{eqn:part-deriv-ac}
\frac{\partial\lambda(t,r)}{\partial t},\ \frac{\partial\lambda(t,r)}{\partial r} \in AC_\BR[0,1],\quad \forall r\in [0,\epsilon_0],
\end{equation}
provided that $\epsilon_0>0$ is sufficiently small. Furthermore,
\begin{equation}\label{eqn:part-deriv-form}
\frac{\partial \lambda(t,r)}{\partial r}=-\Phi^2(t,\lambda(t,r)),\ t\in(0,1),
\end{equation}
where $\Phi(x,\lambda(t,r))$ is  the normalized eigenfunction  of $({\bf E}^{t,r}_q)$ with respect to $\lambda(t,r)$.

\noindent
$(ii)$ $\Phi(x,\lambda(t,r))$ is continuous differentiable with respect to $r\in [0,\epsilon_0]$.
\end{thm}

\begin{prf} 
Since $\lambda(t,r)$ is the first eigenvalue of $({\bf E}^{t,r}_q)$ if and only if
Eq.\,\eqref{eqn:equiv-1st-eigenv} holds, that is $\lambda=\lambda(t,r)$ satisfies the equation
\begin{equation}\label{eqn:F-funct}
F(t,r;\lambda):=r\varphi(t,\lambda)\psi(t,\lambda)+\Delta(\lambda)=0
\end{equation}
for $(t,r)\in (0,1)\times (0,\infty)$. If we choose $\varphi$ and $\psi$ satisfy
\begin{equation}\label{eqn:stand-ini-solut}
 \varphi(0)=0, \ \varphi'(0)=1;\  \psi(1)=0, \ \psi'(1)=-1,
\end{equation}
then the  Wronskian determinant $W[\varphi,\psi](x)$ of $\varphi$ and $\psi$ satisfies
$W[\varphi,\psi](x)\equiv \Delta(\lambda)=-\varphi(1,\lambda)$,
where $\varphi(1,\lambda)$ may be replaced by $\psi(0,\lambda)$ and the proof can be given in the same way, and this $\Delta(\lambda)$ is also regarded as the discriminant function defined as in \eqref{eqn:discrim}.
Consequently, we have
\begin{equation}\label{eqn:F-funct-a}
F(t,r,\lambda)=r\varphi(t,\lambda)\psi(t,\lambda)-\varphi(1,\lambda).
\end{equation}
Set
\begin{equation}\label{eqn:pd-u}
u(x,\lambda)=\frac{\partial\varphi(x,\lambda)}{\partial \lambda}.
\end{equation}
Note that $\lambda(0,r)=\lambda_1$
is the first eigenvalue of the main problem $({\bf E}_q)$ and has simple multiplicity, it holds that
$$
\varphi(1,\lambda(0,r))=\varphi(1,\lambda_1)=0, \ u(1,\lambda_1)\not=0.
$$
Since $\lambda(t,r)$ is continuous on the variable $r\geqslant 0$ and $\lambda(t,0)=\lambda_1$,
we can chose $r$ sufficiently small, say $r\leqslant\epsilon_0$, such that
$$
|\lambda_1-\lambda(t,r)|=|\lambda(t,0)-\lambda(t,r)|
$$
is small enough for $(t,r)\in[0,1]\times [0,\epsilon_0]$
to satisfy $u(1,\lambda(t,r))\not=0$.

On the other hand, we know that $\lambda(t,r)\leqslant\lambda_1$ is continuous for $(t,r)\in[0,1]\times [0,\epsilon_0]$, and
 both of $\varphi(x,\lambda)$ and $\psi(x,\lambda)$ are continuous differentiable on
$(x,\lambda)\in [0,1]\times[\lambda_0,\lambda_1]$, where
$$
\lambda_0=\min\{\lambda(t,r):\ (t,r)\in[0,1]\times[0,\epsilon_0]\}.
$$
Then it follows from Theorem
\ref{thm:cont-1st-eigenv-funct} that $\epsilon_0$ can be chosen small enough such that
$$
r\left|\frac{\partial[\varphi(x,\lambda)\psi(x,\lambda)]}{\partial \lambda}\right|<\min\{|u(1,\lambda(t,r)|: \ t\in[0,1]\}.
$$
As a result, for $(t,r)\in [0,1]\times[0,\epsilon_0]$ and $\lambda\in[\lambda_0,\lambda_1]$, we have
$$
\frac{\partial F(t,r,\lambda)}{\partial \lambda}=r\frac{\partial[\varphi(t,\lambda)\psi(t,\lambda)]}{\partial \lambda}
-u(1,\lambda)\not=0.
$$
It follows from the existence theorem for implicit functions that there exists a unique implicit function $\lambda=\lambda(t,r)$
(Note that the existence has been ensured by the spectral theory) such that
$F(t,r,\lambda(t,r))=0$ and $\lambda(t,r)$ is continuous differentiable with the partial derivative formula as follows:
\begin{equation}\label{eqn:diff-eqn-sys}
\left\{\aligned
&\frac{\partial \lambda(t,r)}{\partial t}\left\{u(1,\lambda)-r\frac{\partial[\varphi(t,\lambda)\psi(t,\lambda)]}{\partial \lambda}\right\}
=r\frac{\partial[\varphi(t,\lambda)\psi(t,\lambda)]}{\partial t},\\
&\frac{\partial \lambda(t,r)}{\partial r}\left\{u(1,\lambda)-r\frac{\partial[\varphi(t,\lambda)\psi(t,\lambda)]}{\partial \lambda}\right\}
=\varphi(t,\lambda)\psi(t,\lambda),
\endaligned\right.
\end{equation}
where $\lambda=\lambda(t,r)$. Since
$$\frac{\partial \varphi(x,\lambda)}{\partial x}, \
\frac{\partial \psi(x,\lambda)}{\partial x},\ \frac{\partial[\varphi(x,\lambda)\psi(x,\lambda)]}{\partial \lambda}\in AC_\BR[0,1],
$$
we can see from \eqref{eqn:diff-eqn-sys} that
$$
\frac{\partial \lambda(t,r)}{\partial t}, \ \frac{\partial \lambda(t,r)}{\partial r}\in AC_\BR[0,1].
$$
This proves the statement \eqref{eqn:part-deriv-ac}.

Similarly to \eqref{eqn:eigenf-perturb-probl}, let  $\Phi$ be the eigenfunction corresponding $\lambda(t,r)$ such that
$$
\Phi(x,\lambda(t,r))=\left\{\aligned
&\varphi(x,\lambda(t,r)),\ &x\in[0,t),\\
&c(t,r)\psi(x,\lambda(t,r)),\  &x\in(t,1],
\endaligned\right.
$$
where $c(t,r)=\frac{\varphi(t,\lambda(t,r))}{\psi(t,\lambda(t,r))}$.
Since both $\varphi(t,\lambda(t,r))$ and $\psi(t,\lambda(t,r))$ are continuous differentiable on $(t,r)$ and
$\psi(t,\lambda(t,r))>0$ for $t\in(0,1)$, we can see that $c(t,r)$ is
continuous differentiable on $(t,r)\in(0,1)\times [0,\epsilon_0]$, and hence $\Phi(x,\lambda(t,r))$ is continuous
differentiable on $(t,r)\in(0,1)\times [0,\epsilon_0]$ and $t\not=x$.
Note that for fixed $x\in(0,1)$, one has
$$
\frac{\partial \Phi(x,\lambda(t,r))}{\partial r}=\left\{\aligned
&\frac{\partial\varphi}{\partial\lambda}\cdot\frac{\partial \lambda(t,r)}{\partial r},\ &t\in(x,1),\\
&\left[\frac{\partial c}{\partial\lambda}\psi+c\frac{\partial\psi}{\partial\lambda}\right]\frac{\partial \lambda(t,r)}{\partial r},\  &t\in(0,x),
\endaligned\right.
$$
where $\varphi=\varphi(x,\lambda(t,r)), \psi=\psi(x,\lambda(t,r))$ and $c=c(t,r)$. One can verify that
$$
\frac{\partial \Phi(x,\lambda(x-0,r))}{\partial r}=\frac{\partial \Phi(x,\lambda(x+0,r))}{\partial r}=
\frac{\partial\varphi}{\partial\lambda}\cdot\frac{\partial \lambda(x,r)}{\partial r},
$$
which means $\Phi(x,\lambda(t,r))$ is continuous differentiable with respect to  $r\in[0,\epsilon_0]$.
This proves the  conclusion $(ii)$.

At last, we prove the formula \eqref{eqn:part-deriv-form}. Set $v=\frac{\partial \Phi}{\partial r}$, then we have
\begin{equation}\label{eqn:Phi}
-\Phi''+[q-r\delta(x-t)]\Phi=\lambda(t,r)w\Phi,\ \Phi(0)=0=\Phi(1),
\end{equation}
and 
\begin{equation}\label{eqn:v}
-v''+[q-r\delta(x-t)]v-\delta(x-t)\Phi=\frac{\partial \lambda(r,t)}{\partial r} w\Phi+\lambda(t,r)w v,\ v(0)=0=v(1).
\end{equation}
Multiplying two sides of the equations \eqref{eqn:Phi} and \eqref{eqn:v} by $v$ and $\Phi$ respectively, and then
integrating their difference over $[0,1]$,  we can verify that
\begin{equation}\label{eqn:pd-form}
\frac{\partial \lambda(t, r)}{\partial r}\int^1_0 w(x)\Phi^2(x,\lambda(t,r)) d x=-\Phi^2(t,\lambda(t,r)).
\end{equation}
If  $\Phi$ is a normalized eigenfunction, then the equation \eqref{eqn:pd-form} is just the desired Eq.\,\eqref{eqn:part-deriv-form}.
The proof is completed.
\end{prf}

\bigskip
\section{Recovery of potential by the first eigenvalue function}\label{sec:potential-recovery}
\medskip

In this section, we reconstruct the potential $q$ by its first eigenvalue function  explicitly. This is also a uniqueness result.

\begin{thm}\label{thm:main-thm}
Let $\lambda(t,r;q)$ be the first eigenvalue function of $({\bf E}_q)$. Then we have that, the potential $q(x)$ is determined by $\lambda(t,r;q)$ uniquely, and can be expressed as
\begin{equation}\label{eqn:reconstruct-formula}
q(x)=\frac{\varphi''_0(x)}{\varphi_0(x)}+\lambda_1 w(x),\
\varphi_0(x)=\sqrt{-\frac{\partial \lambda(x,0;q)}{\partial r}},\ x\in(0,1),
\end{equation}
where $\lambda_1$, which equals to $\lambda(t,0;q)$, is the first eigenvalue of $({\bf E}_q)$. 
\end{thm}

\begin{prf} 
Through this proof, we still simply denote $\lambda(t,r;q)$ by $\lambda(t,r)$ as before.

Letting $r\to0^+$ in the second equation in \eqref{eqn:diff-eqn-sys}, we find that
\begin{equation}\label{eqn:phi-psi-pd}
\frac{\partial \lambda(t,0)}{\partial r}\cdot\frac{\partial \varphi(1,\lambda_1)}{\partial \lambda}
=\varphi(t,\lambda_1)\psi(t,\lambda_1),
\end{equation}
where $\varphi(x,\lambda)$ and $\psi(x,\lambda)$ are defined as in \eqref{eqn:stand-ini-solut}. 
Since $\lambda_1$ is the first eigenvalue of $({\bf E}_q)$, there exists a constant $c$ such that $\varphi(x,\lambda_1)=c\psi(x,\lambda_1)$ which is the eigenfunction corresponding to $\lambda_1$.
This also means that
\begin{equation}\label{eqn:c}
c=\frac{\varphi(x,\lambda_1)}{\psi(x,\lambda_1)},\ x\in(0,1).
\end{equation}
Note that $\varphi(1,\lambda_1)=0=\psi(1,\lambda_1)$ and $\psi'(1,\lambda_1)=-1$. Letting $x\to 1^-$ in the \eqref{eqn:c},
we have from L' H\^opital's rule that
$$c=\frac{\varphi'(1,\lambda_1)}{\psi'(1,\lambda_1)}=-\varphi'(1,\lambda_1)>0,$$
that is
\begin{equation}\label{eqn:phi-psi-lambda1}
\varphi(x,\lambda_1)=-\varphi'(1,\lambda_1)\psi(x,\lambda_1).
\end{equation}

Claim that
\begin{equation}\label{eqn:phi-psi-int}
\int^1_0 w\varphi^2(x)dx=\varphi'(1,\lambda_1)u(1,\lambda_1),
\end{equation}
where $u(x,\lambda)$ is defined as \eqref{eqn:pd-u}.
Indeed, it follows from the definition of $\varphi(x,\lambda)$, i.e.
$$
-\varphi''+q\varphi=\lambda w \varphi,\ \varphi(0)=0, \varphi'(0)=1
$$
that $u$ satisfies
$$
-u''+qu=\lambda w u+w\varphi,\ u(0)=0=u'(0).
$$
Then, by the similar argument in \eqref{eqn:Phi}-\eqref{eqn:pd-form}, the above facts yield that
$$
\varphi''u-u''\varphi=w\varphi^2.
$$
Integrating the above equation over $[0,1]$, we have
$$ \int^1_0 w(x)\varphi(x)^2\,{\rm d}x=\varphi'(1,\lambda)u(1,\lambda)-\varphi(1,\lambda)u'(1,\lambda).$$
So, due to $\varphi(1,\lambda_1)=0$, the above claim holds.

Substituting both of Eq.\,\eqref{eqn:phi-psi-lambda1} and Eq.\,\eqref{eqn:phi-psi-int} into Eq.\,\eqref{eqn:phi-psi-pd}, there has
\begin{equation}\label{44}
\frac{\partial \lambda(t,0)}{\partial r}=-\frac{\varphi^2(t,\lambda_1)}{\int_0^1 w\varphi^2{\rm d}x},\ t\in(0,1).
\end{equation}
Since $\Phi(x,\lambda_1)=\alpha\varphi(x,\lambda_1)$, where $\Phi(x,\lambda_1)$ is a normalized eigenfunction of $({\bf E}_q)$ associated to $\lambda_1$ and $\alpha$ is a constant,
then it follows from Eq.\,\eqref{44} that
\begin{equation}\label{eqn:eigenf-1st-eigenv-funct}
\frac{\partial \lambda(t,0)}{\partial r}=-\Phi^2(t,\lambda_1),\ t\in(0,1).
\end{equation}
Therefore, Eq.\,\eqref{eqn:eigenf-1st-eigenv-funct} gives that
$$
\Phi(t,\lambda_1)=\sqrt{-\frac{\partial \lambda(t,0)}{\partial r}},\ t\in[0,1],
$$
which is exactly the desired $\varphi_0$ in the statement of the present theorem.

Moreover, it has already been proved by Theorem \ref{thm:diff-1st-eigenv-funct} that $\frac{\partial \lambda(t,r)}{\partial t}\in AC_\BR[0,1]$, and hence
$\Phi''(x,\lambda_1)$ exists a.e. on $(0,1)$. Thus, we obtain
$$
q(x)=\frac{\Phi''(x,\lambda_1)}{\Phi(x,\lambda_1)}+\lambda_1 w(x),
$$
which prove that $q(x)$ is able to be reconstructed by $\lambda(t,r)$ uniquely and explictly.
\end{prf}

\begin{rem}\label{rem:approx-recov}
$(i)$ The above theorem also implies the uniqueness:  if $\lambda(t,r;q_1)=\lambda(t,r;q_2)$, then $q_1=q_2$.

\smallskip
\noindent
$(ii)$ Note that $\frac{\partial\lambda(t,0)}{\partial r}$ is a right partial derivative at $r=0$.
In fact, if $\lambda(t,r)$ is the first eigenvalue function of $({\bf E}_q)$ with an unknown potential $q$, then, to recover $q$, we only need to know the values of
$\lambda(t,r_n)$ for $r_n\to 0^+$ as $n\to\infty$ and $t\in (0,1)$, since one has
$$\frac{\partial \lambda(t,0)}{\partial r}=\lim\limits_{n\to\infty}\frac{\lambda(t,r_n)-\lambda(t,0)}{r_n}.$$
\end{rem}

By now, Theorem~\ref{thm:main-thm}, together with Remark~\ref{rem:approx-recov}$(i)$, has already provided  the uniqueness and reconstruction formula for inverse S-L problem proposed at the beginning of our paper. 
But, there still remain many unsolved questions, such as existence of the first eigenvalue function, stability of reconstructing the potential, and so on (cf. \cite{FY2008, FWW2016}).  With the end of this paper, we are going to give a few discussions on these topics.

On one hand, from Theorems \ref{thm: equiv-1st-eigenv}--\ref{thm:main-thm}, we can see that the first eigenvalue function $\lambda(t,r)$ associated to $({\bf E}_q)$ needs to satisfy the following conditions:

{\it
\smallskip
\noindent
$(i)$ $$\frac{\partial \lambda(t,r)}{\partial r}, \ \frac{\partial \lambda(t,r)}{\partial t}\in AC_\BR[0,1],$$
for a sufficiently small $\epsilon_0>0$ such that $r\in [0,\epsilon_0]$;

\smallskip
\noindent
$(ii)$ $$\frac{\partial \lambda(t,r)}{\partial r}<0, \quad\int^1_0 \frac{\partial \lambda(t,0)}{\partial r}dt=-1;$$

\smallskip
\noindent
$(iii)$ $$ \lambda(t,0)=\lambda(0,r)=\lambda(1,r):=\lambda_1\geqslant\lambda(t,r);$$

\smallskip
\noindent
$(iv)$ Define
\begin{equation}\label{322qq}
q_0(x)=\frac{\varphi''_0(x)}{\varphi_0(x)}+\lambda_1 w(x),\quad
\varphi_0(x)=\sqrt{-\frac{\partial \lambda(x,0)}{\partial r}},
\end{equation}
where $\varphi(x,\lambda)$ and $\psi(x,\lambda)$ are the solutions of
$-y''+q_0 y=\lambda wy$ such that
$$
\varphi(0,\lambda)=0,\ \varphi'(0,\lambda)=1;\ \psi(1,\lambda)=0,\ \psi'(1,\lambda)=-1.
$$
Then $\lambda(t,r)$ should be the unique root of
$$
r\varphi(t,\lambda)\psi(t,\lambda)-\varphi(1,\lambda)=0.
$$
Furthermore, the corresponding potential
is the  $q_0$ defined in \eqref{322qq}. 
}

\bigskip

For example, 
$$
\lambda(t,r)=\lambda_1-r\sin^2(\pi t),\ \lambda_1\in (-\infty,\infty),
$$
then it can not be the first eigenvalue function of any problem $({\bf E}_q)$, since the condition $(ii)$
is not satisfied. 

The following is a motivating question of the discussion above.

\medskip
\noindent
{\bf Question 1.} \emph{What assumptions should $\lambda(r,t)$ satisfy such that it is the first eigenvalue function associated to some problem $({\bf E}_q)$?}
\medskip

On the other hand, we look at  an example as below.

\begin{eg}\label{ex:0-potential}
Suppose that $\lambda(t,r)=\rho^2(t,r)$ is the first eigenvalue function of $({\bf E}_q)$ with $q(x)\equiv 0$ as well as $w\equiv 1$. Then, by Theorem~\ref{thm: equiv-1st-eigenv}, the function $\rho(t,r)$ must satisfy the following non-linear equation
\begin{equation}\label{eqn:1st-eigenv-funct-q=0}
\rho\sin\rho=r\sin(\rho t)\sin(\rho(1-t)),\ t\in(0,1)
\end{equation}
with the conditions that $\rho(0,r)=\rho(1,r)=\rho(t,0)=\pi$, or equivalently,
\begin{equation}\label{eqn:1st-eigenv-funct-q=0-a}
\rho[\cot(\rho t)+\cot(\rho(1-t))]=r,\ t\in(0,1),\ \rho(t,0)=\pi,
\end{equation}
where $r\in[0,\epsilon_0]$ with  sufficient small enough $\epsilon_0>0$, and $\rho=\rho(t,r)\in (0,\pi)$ for $t\in(0,1)$.  In fact, one can verify
\begin{equation}\label{eqn:1st-eigenv-funct-q=0-b}
\frac{\partial \rho(t, r)}{\partial r}
\left(\cot(\rho t)+\cot(\rho(1-t))-\rho\left[\frac{t}{\sin^2(\rho t)}+\frac{1-t}{\sin^2(\rho(1-t))}\right]\right)=1.
\end{equation}
Let $r=0$ in the above equation, and note that $\rho(t,0)=\pi$. Then we find that for $t\in(0,1)$,
$$
\pi\frac{\partial \rho(t, 0)}{\partial r}
\left[\frac{t}{\sin^2(\pi t)}+\frac{1-t}{\sin^2(\pi(1-t))}\right]=-1.
$$
This gives that
\begin{equation}\label{eqn:eigenf-1st-eigenv-funct-q=0}
\frac{\partial \lambda(t, 0)}{\partial r}=2\rho(t,0)\frac{\partial \rho(t, 0)}{\partial r}=-2\sin^2(\pi t):=-\Phi^2(t,\lambda_1),
\end{equation}
where $\lambda_1=\pi^2$. 
It is obvious that $\Phi(t,\lambda_1)=\sqrt{2}\sin(\pi x)$ is a normalized eigenfunction of this eigenvalue problem with respect to the first eigenvalue $\lambda_1$.
So, Eq.\,\eqref{eqn:eigenf-1st-eigenv-funct-q=0} is just Eq.\,\eqref{eqn:eigenf-1st-eigenv-funct} for the present example, and it can be seen from both of Eq.\,\eqref{eqn:eigenf-1st-eigenv-funct-q=0} and Theorem~\ref{thm:main-thm} that the zero potential is exactly recovered by $\lambda(t,r)$. 
\end{eg}

Suppose that $\lambda(t,r,q_1)$ and $\lambda(t,r,q_2)$ are the first eigenvalue functions of $({\bf E}_{q_1})$ and $({\bf E}_{q_2})$, respectively.
From Theorem \ref{thm:main-thm} and Example~\ref{ex:0-potential}, we can see that, if
$$\lambda(t,0,q_1)=\lambda(t,0,q_2):=\lambda_1, \ \frac{\partial \lambda(t,0;q_1)}{\partial r}=\frac{\partial \lambda(t,0;q_2)}{\partial r}, \ t\in(0,1),$$
then  $q_1=q_2$, and hence we have $\lambda(t,r,q_1)=\lambda(t,r,q_2)$. 
Note that $$\lambda(0,r,q_j)=\lambda(1,r,q_j)=\lambda(t,0,q_j),\ j=1,2.$$ 
The above facts mean that a first eigenvalue function $\lambda(t,r)$ for certain SLP can be determined by
\begin{equation}\label{eqn:ini-value-1st-eigenv-funct}
\lambda(0,r)=\lambda(1,r)=\lambda(t,0)\equiv\text{constant} \text{  and  }\frac{\partial \lambda(t,0)}{\partial r}<0.
\end{equation}

However, the following example implies that, it is probable that a function $f(t,r)$, satisfying the conditions \eqref{eqn:ini-value-1st-eigenv-funct}, is not the first eigenvalue of any problem $({\bf E}_q)$.

\begin{eg}\label{ex:counterex}
$$
\lambda(t,r):=\pi^2-2r\sin^2(\pi t)
$$
is not the first eigenvalue function of  $({\bf E}_q)$ with $w\equiv 1$. 

In fact, suppose that  $\lambda(t,r)$ is the first eigenvalue function. Then,  by Eq.\,\eqref{eqn:diff-eqn-sys}, it can be seen that, the function $\lambda(t,r)$ should satisfy the following equation
\begin{equation}\label{eqn:F-phi-psi}
\frac{\partial \lambda(t,r)}{\partial r}\cdot\frac{\partial F}{\partial \lambda}
=\varphi(t,\lambda)\psi(t,\lambda),
\end{equation}
where $F$ is defined by Eq.\,\eqref{eqn:F-funct}, and
$$
\frac{\partial F}{\partial \lambda}=u(1,\lambda)-r\frac{\partial[\varphi(t,\lambda)\psi(t,\lambda)]}{\partial \lambda},\
u(1,\lambda)=\frac{\partial \varphi(1,\lambda)}{\partial \lambda}.
$$
Since
$$
\frac{\partial \lambda(t,r)}{\partial r}=-2\sin^2(\pi t),
$$
we can see from Example~\ref{ex:0-potential} that if  such a function $\lambda(t,r)$ is the first eigenvalue function of certain SLP $({\bf E}_q)$,
then the corresponding potential $q$ must be $0$. As a result,  we can choose $\varphi$ and $\psi$ in Eq.\,\eqref{eqn:F-phi-psi} as follows:
\begin{equation}\label{eqn:stand-int-solut-q=0}
\varphi(t,\lambda)=\frac{1}{\sqrt{\lambda}}\sin\sqrt{\lambda}t,\
\psi(t,\lambda)=\frac{1}{\sqrt{\lambda}}\sin\sqrt{\lambda}(1-t).
\end{equation}
Taking $t=\frac{1}{2}$, by \eqref{eqn:stand-int-solut-q=0}, we find that Eq.\,\eqref{eqn:F-phi-psi} does not hold.  
This means that $\lambda(t,r)=\pi^2-2r\sin^2(\pi t)$ never becomes a first eigenvalue function of any problem $({\bf E}_q)$ with $w\equiv 1$.
\end{eg}

Hence, combining the equations \eqref{eqn:1st-eigenv-funct-q=0-b}--\eqref{eqn:F-phi-psi}, we reasonably guess that $\lambda(t,r)$ maybe satisfy certain initial-boundary values problem of second order partial differential equation
involving $(t,r)$ such as 
\begin{equation}\label{eqn:Q2}
\left\{\aligned
& P\left(\frac{\partial^2 \lambda(t,r)}{\partial t^2},\frac{\partial^2 \lambda(t,r)}{\partial t\partial r},
\frac{\partial \lambda(t,r)}{\partial t},\frac{\partial \lambda(t,r)}{\partial r},\lambda(t,r)\right)=0,\\
& \lambda(0,r)=\lambda(1,r)=\lambda(t,0):=\lambda_1,\ \frac{\partial \lambda(t,0)}{\partial r}=f(t)<0.
\endaligned\right.
\end{equation}

The following is a natural question arising from the examples above.

\medskip
\noindent
{\bf Question 2.}
\emph{Can $\lambda(t,r)$ be characterized by a solution of certain initial-boundary value problem of second order partial differential equation involving $(t,r)$? }
\medskip

As long as either of the questions above can be solved, then the existence of the first eigenvalue function is gained.

\medskip

Here is another question related to Remark~\ref{rem:approx-recov}.

\medskip
\noindent
{\bf Question 3.}  \emph{If only the values of $\lambda(t_m,r_n)$  are observed such that
$$r_n>0,\ n\in\BZ^+,\ \lim_{n\rightarrow +\infty}r_n=0; \ \{t_m:\ 1\leqslant m\leqslant N,\ N\in\BZ^+\} \text{ is a finite subset of } [0,1].$$
Then how can we recover $q$ approximately? Can we estimate the errors? This question involves the stability of recovering the unknown potential of our problem.}

\bigskip

\appendix \section {The proof of Theorem~\ref{thm:cont-eigenv}}\label{sec:app}
\medskip

In this appendix, we will prove Theorem~\ref{thm:cont-eigenv}.  Before the proof, as preparation, we need to introduce some concepts as well as known facts, and also prove several auxiliary results that are crucial for our proof. 
As the generalizations of the classical theorems, Theorem~\ref{thm:eigenv-simply} and Theorem~\ref{thm:zero-eigenf} prove the geometric multiplicity of eigenvalues and  oscillation of  eigenfunctions for $({\bf E}_\mu)$  with positive integrable weights, respectively.

\begin{thm}\label{prop:montel}\cite[Theorem 3.3]{SS2003}(Montel's Theorem) Suppose $\CF$ is a family of holomorphic functions on $\Omega$ that is uniformly bounded on compact subsets of a domain $\Omega\subset\BC$, where $\Omega$ is an open subset of $\BC$.
Then one has

\smallskip
\noindent
$(i)$ $\CF$ is equicontinuous on every compact subset of $\Omega$;

\smallskip
\noindent
$(ii)$ every sequence in $\CF$ has a subsequence that converges uniformly on every compact subset of $\Omega$ (the limit need not be in $\CF$).
\end{thm}
Here, a domain in $\BC$ is an open subset of complex plane. The family $\CF$ is called to be {\bf uniformly bounded on compact subsets} (resp. {\bf on a bounded subset}) of $\Omega$ if for each compact set  (resp. a bounded set) $K\subset\Omega$, there exists $B>0$, such that $\abs{f(z)}\leqslant B$ for all $z\in K$ and $f\in\CF$.
Also, the family $\CF$ is called to be {\bf equicontinuous on a compact set} $K$ if for every $\varepsilon>0$ there exists $\delta>0$ such that whenever $z_1, z_2\in K$ and $\abs{z_1-z_2}<\delta$, then one has $\abs{f(z_1)-f(z_2)}<\varepsilon$ for all $f\in\CF$. See \cite[Page 225]{SS2003} and so on.

As an application, we have the following theorem about the zeros of a holomorphic function.
\begin{lem}\cite[Lemma 2.2]{XGL2018}\label{lem:zero-holom-funct}
Let $F(\lambda)\not\equiv 0$ and $\{F_n(\lambda)\}_{n\in\BZ^+}$ be holomorphic functions on $\mathbb{C}$.
Suppose that $\{\abs{F_n(\lambda)}\}_{n\in\BZ^+}$ is uniformly bounded on compact subsets of $\BC$, and $\lim_{n\rightarrow +\infty}F_n(\lambda)=F(\lambda)$ for any $\lambda\in\BC$.
Denote by $\Sigma_n$ and $\Sigma$ be the zero point sets of $F_n$ and $F$, respectively.
Then one has

\smallskip
\noindent
$(1)$ $\Sigma=\{\lambda:\ \exists \lambda_n\in\Sigma_n,\ s.t.\ \lim_{n\rightarrow +\infty}\lambda_n=\lambda\}.$

\smallskip
\noindent
$(2)$ If $\lambda_0\in\Sigma$, and for every $n\in\BZ^+$, there exists distinct $\lambda_{n,1}, \lambda_{n,2}\in \Sigma_n$ such that $\lim_{n\to \infty} \lambda_{n,j}=\lambda_0\ (j=1,2)$, then one has $F'(\lambda_0)=0$.
\end{lem}

In the rest of this part, unless otherwise specified, we simply denote by $\mu(\lambda)$ the measure $\mu(x,\lambda)\in \CM^\BK_0$ on $[0,1]$, where $x\in[0,1]$ is the argument and $\lambda\in\BC$ is a parameter.

\begin{thm}\label{thm: unifm-cont-paramtr}
Let $\{\mu_n(\lambda)\}_{n\in\BZ^+,\lambda\in\BC}$ and $\{\mu_0(\lambda)\}_{\lambda\in\BC}$ be contained in $\CM^\BK_{0}$ satisfying
$$\mu_n(\lambda)\stackrel{w^{\ast}}{\longrightarrow} \mu_0(\lambda),\quad \text{as }n\rightarrow +\infty, \quad  \forall \lambda\in\BC,$$ and suppose that the sequence $\{\norm{\mu_n(\lambda)}_{\bf V}\}_{n\in\BZ^+}$ is uniformly bounded on any bounded domain $\Omega\subset \mathbb{C}$.
Denote by $y(x,\mu_n(\lambda))$ the solution of the problem \eqref{eqn:mde} with $\mu=\mu_n(\lambda).$

For any fixed $x\in[0,1]$, if $y(x,\mu_n(\lambda))$ and $y^{\bullet}(1,\mu_n(\lambda))$ are holomorphic functions with respect to $\lambda \in \mathbb{C}$, then
$y(x,\mu_n(\lambda))$ and $y^{\bullet}(1,\mu_n(\lambda))$ respectively converge uniformly to $y(x,\mu_0(\lambda))$ and $y^{\bullet}(1,\mu_0(\lambda))$ with respect to $\lambda$ on any compact subset $K$ of $\Omega$, as $n$ tends to $+\infty$.
\end{thm}

\begin{prf}
Following the similar lines of argument as the proof of \cite[Lemma 3.10]{MZ2013}, one has
\begin{equation}\label{eqn:y-bound}
|y(x,\mu_n(\lambda))|\leqslant(|y_0|+|z_0|)e^{\sup\{\|\mu_n(\lambda)\|_{\bf V}:\ n\geqslant1,\ \lambda \in \Omega\}},\ \forall n\geq1,\ x\in[0,1].
\end{equation}
Since the sequence $\{\norm{\mu_n(\lambda)}_{\bf V}\}_{n\in\BZ^+}$ is uniformly bounded on $\Omega$, by the inequality \eqref{eqn:y-bound}, we can know that the sequence $\{\norm{y(\cdot,\mu_n(\lambda))}_\infty\}_{n\in\BZ^+}$ is also uniformly bounded on $\Omega$.
Since every $y(x,\mu_n(\lambda))$ is a holomorphic function with respect to $\lambda\in\Omega$ for any fixed $x$ and all $n\in\BZ^+$, we can see from Theorem~\ref{prop:montel}  that for any fixed $x$, every subsequence in $\{y(x,\mu_n(\lambda))\}_{n\in\BZ^+}$ has a convergent subsequence  that converges uniformly on any compact subset $K$ of $\Omega$.
From \cite[Theorem 1.1]{MZ2013}, it follows that
$y(x,\mu_n(\lambda))$ converges uniformly to $y(x,\mu_0(\lambda))$ with respect to $\lambda$ on any compact subset $K$ of $\Omega$ as $n$ tends to infinity.
Similarly, we also can prove that $y^{\bullet}(1,\mu_n(\lambda))$ converges uniformly to $y^{\bullet}(1,\mu_0(\lambda))$ on any compact subset $K$ of $\Omega$.
\end{prf}

Recall the MDE eigenvalue problem $({\bf E}_\mu)$ defined as in \eqref{eqn:mu-probl-mde}, which is equivalent to the following
system for $(y(x),y^{\bullet}(x))$
\begin{equation}\label{eqn:mu-equiv-system}
{\rm d}y(x)=y^{\bullet}(x){\rm d}x,\
{\rm d}y^{\bullet}(x)=y(x){\rm d}\mu(x)-\lambda wy{\rm d}x,
\end{equation}
together with boundary condition
\begin{equation}
AY(0)+BY(1)=0,\
Y=\begin{pmatrix}
y\\
y^{\bullet}
\end{pmatrix},\
A=
\begin{pmatrix}
1 & 0\\
0 & 0
\end{pmatrix},\
B=
\begin{pmatrix}
0 & 0\\
1 & 0
\end{pmatrix}.
\end{equation}
Denote by $y(x,\lambda;y_0,z_0)$ the unique solution of the system of equations \eqref{eqn:mu-equiv-system} with the initial values $(y(x),y^{\bullet}(x))=(y_0,z_0)$.

Let
\begin{equation}
\Phi(x,\lambda):=\begin{pmatrix}
\varphi(x,\lambda) & \psi(x,\lambda)\\
\varphi^{\bullet}(x,\lambda) & \psi^{\bullet}(x,\lambda)
\end{pmatrix}
\end{equation}
be the primary fundamental matrix of \eqref{eqn:mu-equiv-system}, where $\varphi(x,\lambda):= y(x,\lambda;1,0),\ \psi(x,\lambda):= y(x,\lambda;0,1)$.
Define the {\bf discriminant function} for $({\bf E}_\mu)$ by
\begin{equation}\label{eqn:discrim}
F(\lambda)=det[A+B\Phi(1,\lambda)].
\end{equation}

\begin{rem}\label{rem:eigenv-F-zero}
Here $F(\lambda)=\psi(1,\lambda)$, and note that $\lambda$ is an eigenvalue of the problem $({\bf E}_\mu)$ if and only if $F(\lambda)=0$. 

In some other literatures, the discriminant function is also called characteristic function, which could be confused with the function $\chi_U$ defined at the beginning of Section~\ref{sec:prelim}. For more details, refer to \cite[Page~45]{Zettl2005} and \cite{MZ2013}.
\end{rem}

\begin{rem}\label{rem:discrim-funct}
It is well known that for the problem $({\bf E}_\mu)$ with a.c. potential $\mu$, i.e., a classical Dirichlet eigenvalue problem,  its discriminant function is holomorphic. Hence, from both of Remark \ref{rem:L1-approx} and Theorem \ref{thm: unifm-cont-paramtr}, it follows that the discriminant function $F(\lambda)(=\psi(1,\lambda))$ for $({\bf E}_\mu)$ is also a holomorphic function, since the limit function of a uniformly convergent sequence of holomorphic functions is still a holomorphic function on any bounded domain.

\end{rem}

The following theorem generalizes the multiplicity theorem of eigenvalues for the classical S-L problems.

\begin{thm}\label{thm:eigenv-simply}
The measure eigenvalue problem $({\bf E}_\mu)$ exists countable algebraically simple eigenvalues
$$-\infty<\lambda_1<\cdots<\lambda_n\to\infty,\ n\to\infty.$$
\end{thm}
\begin{prf} 
From Remark~\ref{rem:L1-approx}$(ii)$, for the given $\mu$ in $({\bf E}_\mu)$, we can choose $\{q_n\}_{n\in\BZ^+}\subset L^1[0,1]$ such that
$\mu_{q_n} \stackrel{w^{\ast}}{\longrightarrow} \mu$, and thus $\{\norm{q_n}_{L^1}\}_{n\in\BZ^+}$ is bounded.

Let $F(\lambda)$ and $F_n(\lambda)$ be the discriminant functions for $({\bf E}_\mu)$ and $({\bf E}_{\mu_{q_n}})$, respectively. So Remark~\ref{rem:eigenv-F-zero} tells us that $F(\lambda)=\psi(1,\mu-\lambda\mu_w)$ and $F_n(\lambda)=\psi(1,\mu_{q_n}-\lambda\mu_w)$ are the solutions of the initial value problems \eqref{eqn:mde} with the measure $\mu-\lambda\mu_w$ and $\mu_{q_n}-\lambda\mu_w$ under the initial values $(0,1)$, respectively.
Then, by Theorem \ref{thm: unifm-cont-paramtr}, we have that $F_n(\lambda)$ converges uniformly to $F(\lambda)$ with respect to $\lambda$ on any compact subset $K$ of bounded domain $\Omega\subset \mathbb{C}$ as $n$ going to $+\infty$.

Let $\lambda_m(q_n)(m\in \BZ^+)$ be the $m$-th eigenvalue of  $({\bf E}_{\mu_{q_n}})$ which is a classical S-L problem. 
According to \cite[Theorem VII.3.6]{Kato1980}, which gives the boundedness of eigenvalues of SLPs with integrable bounded potentials,  we know that $\{\lambda_m(q_n)\}$ is bounded for any $m$.
Thus $\{\lambda_m(q_n)\}_{n\in\BZ^+}$ has a convergent subsequence for any fixed $m$.
 We may as well assume that $\{\lambda_m(q_n)\}_{n\in\BZ^+}$ is convergent, and the limit is denoted by $\alpha_m$.
Since $F_n(\lambda)$ converges uniformly to $F(\lambda)$,  we know that $F'_n(\lambda)$ is bounded uniformly.
Since
$$\abs{F_n(\alpha_m)-F_n(\lambda_m(q_n))}\leqslant \abs{F'_n(\xi)}\abs{|\alpha_m-\lambda_m(q_n)}\rightarrow0, \text{ as }n\rightarrow +\infty,$$
and $F_n(\lambda_m(q_n))=0$ (Remark~\ref{rem:eigenv-F-zero}), it holds that  $F_n(\alpha_m)=0$. 
Immediately, it follows from the uniform convergence of $F_n$ that $F(\alpha_m)=0$, by Remark~\ref{rem:eigenv-F-zero} again, means that $\alpha_m$ is the eigenvalue of the problem $({\bf E}_\mu)$. 
By Lemma~\ref{lem:zero-holom-funct}$(1)$, we can know that $\{\alpha_m\}_{m\in\BZ^+}$ is exactly the eigenvalue sequence $\{\lambda_m\}_{m\in\BZ^+}$ of $({\bf E}_\mu)$.
Hence, by the boundedness below and infinite increasing of $\{\lambda_m(q_n)\}$ for any $n$, it can be seen that $\{\lambda_m\}_{n\in\BZ^+}$ is an infinitely countable sequence, which is bounded below and tends to $+\infty$ as $n$ goes to $+\infty$.

On the other hand, similarly as the proof of the \cite[Theorem III.1.2]{Min1983}, it can be proved that all eigenvalues of the problem \eqref{eqn:mu-probl-mde} are real numbers, which are the zeros of $F(\lambda)$ by Remark~\ref{rem:eigenv-F-zero}. 
Since zeros of holomorphic function are isolated, one has that $F(\lambda)$ only has countable zeros.
Up to now, we have proved that $({\bf E}_\mu)$ has an ascending sequence of infinitely countable real eigenvalues, which is bounded below.

It remains to prove that each eigenvalue is algebraically simple (automatically, geometrically simple).
From Remark \ref{rem:discrim-funct}, there exists
$\psi_\lambda=\frac{\partial \psi}{\partial \lambda}$ such that
\begin{equation}\label{eqn:psi-lambda}
\left\{\aligned
&-{\rm d}\psi^{\bullet}_\lambda(x)+\psi_\lambda(x){\rm d}\mu(x)=\lambda w \psi_\lambda+w \psi,\ x\in [0,1],\\
&\psi_\lambda(0)=0=\psi^{\bullet}_\lambda(0).
\endaligned\right.
\end{equation}
Multiplying two sides of the equation \eqref{eqn:psi-lambda} and \eqref{eqn:mu-probl-mde} by $\psi$ and $\psi_\lambda$ respectively, and then
integrating their difference over $[0,1]$,
we have
\begin{equation}\label{eqn:psi-lambda-2}
\int_0^1 w|\psi|^2{\rm d}x=
\psi_\lambda\psi^{\bullet}(1)-\psi\psi^{\bullet}_\lambda(1).
\end{equation}
Since $F(\lambda)=\psi(1,\lambda)=0$ and $\int_0^1 w|\psi|^2{\rm d}x\neq 0$,  we obtain by Eq.\,\eqref{eqn:psi-lambda-2} that $F'(\lambda)=\psi_\lambda\neq 0$, which means each eigenvalue is algebraically simple.
\end{prf}

\begin{lem}\label{lem:zero-cont-funct}
Let $y(x)$ be the eigenfunction of the problem \eqref{eqn:mu-probl-mde}, and $\{y_n(x)\}_{n\in\BZ^+}$ be a sequence of continuous functions such that converges uniformly to $y(x)$ on $[0,1]$.
If every $y_n(x)$ has exactly $m\in \mathbb{N}$ zeros on $(0,1)$, then $y(x)$ also has only  $m$ zeros on $(0,1)$.
\end{lem}
\begin{prf} Let $x_{1n}, x_{2n}, \cdots, x_{mn}$ be $m$ zeros of $y_n(x)$ on the open interval $(0,1)$.  We may as well assume that $\lim_{n\to \infty} x_{kn}=x_{k},\ 1\leq k\leqslant m.$
From the properties of uniformly convergent continuous functions, it follows that $y(x_k)=0$.

Set $I_k=(x_k,x_{k+1}),\ 0\leq k\leqslant m,$ where $x_0=0,\ x_{m+1}=1$.
It will be enough to prove that $y(x)$ does not have sign-change on $I_k,\ 0\leqslant k\leqslant m$, or equivalently, to show that, for any sufficiently small $\delta>0$, one always has $y(x)>0$ (or $y(x)<0$) on $I_{k\delta}=(x_k+\delta, x_{k+1}-\delta),\ 0\leqslant k\leqslant m$.

We may as well assume that $y_n(x)>0$ for all $x\in I_{k\delta}$, then $y(x)\geqslant 0$ for all $x\in I_{k\delta}$.
If there exists $x_0 \in I_{k\delta}$ such that $y(x_0)=0$, from the definition \eqref{eqn:solution-mde} of solutions of a MDE, one can verify that
$$y'(x_{0}^{+})-y'(x_0^{-})=\lim_{\epsilon\rightarrow 0^+}\int^{x_0+\epsilon}_{x_0-\epsilon} y(x)\,{\rm d} \{\mu(x)-\lambda \mu_w(x)\}=0,$$
which means that $y'(x_0)=0.$
Then $y(x)$ is the zero solution of the initial value problem
$$-{\rm d}y^{\bullet}(x)+y(x)\,{\rm d}\mu(x)=\lambda y(x) w(x)\,{\rm d}x,\ x\in [0,1],\ y(x_0)=0=y^{\bullet}(x_0),$$
which is a contradiction.\end{prf}

The number of zeros of eigenfunctions of the problem $({\bf E}_\mu)$ is computed out in the next theorem.

\begin{thm}\label{thm:zero-eigenf}
For any $m\in \BZ^+$, the $m$-th eigenfunction $y_m(x)$ of the problem $({\bf E}_\mu)$ has exactly $m-1$ zeros  on $(0,1)$.
\end{thm}

\begin{prf}
From Remark~\ref{rem:L1-approx}$(ii)$, we know that, for the given $\mu$ in $({\bf E}_\mu)$, there exists $q_n\in L^1[0,1]$ such that
$\mu_{q_n} \stackrel{w^{\ast}}{\longrightarrow} \mu$.
Set $F_n(\lambda)$ and $F(\lambda)$ as those in the proof of Theorem~\ref{thm:eigenv-simply}, so they satisfy the conditions of Lemma \ref{lem:zero-holom-funct}. 
Let  $\lambda_m(\mu_{q_n})$ and $\lambda_m(\mu)$ be the $m$-th eigenvalue of $({\bf E}_{\mu_{q_n}})$ and $({\bf E}_\mu)$, respectively.
So  $\{\lambda_m(\mu_{q_n})\}_{m\in\BZ^+}$ and $\{\lambda_m(\mu)\}_{m\in\BZ^+}$ are the sets of zeros of $F_n(\lambda)$ and $F(\lambda)$, respectively.
We firstly prove $\lim_{n\to \infty} \lambda_m(\mu_{q_n})=\lambda_m(\mu)$ for any $m\geqslant 1$.

For $m=1$, note that $\{\lambda_1(\mu_{q_n})\}_{n\in\BZ^+}$ is bounded according to \cite[Theorem VII.3.6]{Kato1980}. Suppose there exist a subsequence of $\{\lambda_1(\mu_{q_n})\}_{n\in\BZ^+}$ that does not converge to $\lambda_1(\mu)$. We may as well assume that $\lim_{n\to \infty} \lambda_1(\mu_{q_n})=\alpha_1\neq\lambda_1(\mu)$.
Then there exists $m>1$ such that $\lambda_m(\mu)=\alpha_1>\lambda_1(\mu)$ by Lemma \ref{lem:zero-holom-funct}$(i)$ and Theorem~\ref{thm:eigenv-simply}.
Since $\lambda_m(\mu_{q_n})>\lambda_1(\mu_{q_n})$ for any $m>1$, we have that any accumulation point of $\{\lambda_m(\mu_{q_n})\}_{n\in\BZ^+}$ for any $m>1$ is strictly bigger than $\lambda_1(\mu)$, and then it turns out that $\lambda_1(\mu)$ is not the eigenvalue of the problem $({\bf E}_\mu)$ by Lemma \ref{lem:zero-holom-funct}$(i)$, which is a contradiction.

The similar argument is applied to the case $m=2$. Suppose $\lim_{n\to \infty} \lambda_2(\mu_{q_n})=\beta_2\neq\lambda_2(\mu)$.
According to Lemma \ref{lem:zero-holom-funct}$(i)$, we know that $\beta_2=\lambda_1(\mu)$ or there exists $m>3$ such that $\beta_2=\lambda_m(\mu)$. The first case contradicts to  Theorem \ref{thm:eigenv-simply} due to Lemma \ref{lem:zero-holom-funct}$(ii)$.
For the second case, by the similar procedures in the above case $m=1$, we also can get that $\lambda_2(\mu)$ is not an eigenvalue of the problem \eqref{eqn:mu-probl-mde}, which is also a contradiction.
Then, by induction, we obtain that $\lim\limits_{n\to \infty} \lambda_m(\mu_{q_n})=\lambda_m(\mu)$ for any $m\in\BZ^+$.

Now, if $y_m(x;\mu_{q_n})$ is the $m$-th eigenfunction of the problem $({\bf E}_{\mu_{q_n}})$ and has the same initial condition at $x=0$ as that of $y_m$, 
then $y_m(x;\mu_{q_n})$ uniformly converges to $y_m$ with respect to $x\in[0,1]$ by \cite[Theorem 1.1]{MZ2013}.
Since for any $m\geqslant 1$,  the eigenfunction $y_m(x;\mu_{q_n})$ has $m-1$ zeros on $(0,1)$ from \cite[Theorem 4.3.1]{Zettl2005}, the eigenfunction $y_m$ also has $m-1$ zeros on $(0,1)$ by Lemma \ref{lem:zero-cont-funct}.
\end{prf}

Now, it is time for us to give the proof of Theorem~\ref{thm:cont-eigenv}.

\medskip
\noindent
{\bf The proof of Theorem~\ref{thm:cont-eigenv}.} 
Suppose that a sequence $\{\mu_n\}_{n\in\BZ^+}\subset\CM^\BR_0$ is weakly$^*$ convergent to $\mu_0\in\CM^\BR_0$. 
Lemma \ref{lem:zero-holom-funct} says that, for any $m\in \BZ^+$, there exists $k\in \BZ^+$ such that $$\lim_{n\to \infty} \lambda_k(\mu_{n})=\lambda_m(\mu_0).$$ 
Then the corresponding eigenfunction $y_k(\mu_n)$ converges uniformly to $y_m(\mu_0)$ with respect to $x\in[0,1]$ by \cite[Theorem 1.1]{MZ2013}.
According to Theorem \ref{thm:zero-eigenf}, we can find that $y_k(\mu_n)$ and  $y_m(\mu_0)$ have $k-1$ and $m-1$ zeros on $(0,1)$ respectively.
Then, it follows from Lemma \ref{lem:zero-cont-funct} that $m$ equals to $k$, and so the proof is done.

\bigskip

\bigskip

\section*{Acknowledgement}
\smallskip

This research was partially supported by the NSF of China [Grant numbers 12271299, 12071254, and 11701327].

\end{document}